\documentclass[12pt]{article}
\usepackage[mathscr]{eucal}
\usepackage{amsmath}
\usepackage{amssymb}

\newtheorem{thm}{Theorem}
\newtheorem{rmk}{Remark}
\newtheorem{lem}{Lemma}

\begin{document}

\noindent \textbf{\large AVERAGE SAMPLING RESTORATION OF HARMONIZABLE
PROCESSES\footnotetext{Submitted: \today}}
\vskip 5mm

\noindent {\it Short title:} \textbf{AVERAGE SAMPLING OF HARMONIZABLE PROCESSES}

\vskip 1cm
\noindent \textbf{\large Andriy Olenko$^a$ and Tibor Pog\'any$^b$\footnote{Address correspondence to Tibor Pog\'any, Faculty of
Maritime Studies, University of Rijeka, Studentska 2, HR-51000 Rijeka, Croatia; E-mail: poganj@pfri.hr}}

\vskip 5mm
\noindent {$^a$ Department of Mathematics and Statistics, La Trobe University, Melbourne, Australia}

\noindent {$^b$ Faculty of Maritime Studies, University of Rijeka, Rijeka, Croatia}

\

 This is an Author's Accepted Manuscript of an article published in the
Communications in Statistics - Theory and Methods, Vol. 40, Issue 19-20,  2011,  3587-3598.  [copyright Taylor \& Francis], available online at: http://www.tandfonline.com/ [DOI: 10.1080/03610926.2011.581180]

\vskip 3cm

\vskip 3mm

\vskip 6mm

\noindent {\bf ABSTRACT}

\noindent \textit{The harmonizable Piranashvili--type stochastic processes are approximated by finite time shifted average sampling sums. Explicit
truncation error upper bounds are established. Va\-rious corollaries
and special cases are discussed.}
\vskip 4mm

\newpage


\section{Introduction}
Recovering a stochastic signal from discrete samples or assessing
the information lost in the sampling process are the  fundamental
problems in sampling and interpolation theory. An essential question of the theory of approximation of random functions is the problem of
determining the classes of approximant.

We start with the Paley-Wiener class $PW^2_{\omega}$ of all non-random complex-valued $L^2(\mathbb R)$-functions whose Fourier
spectrum is bandlimited to $[-\omega,\omega)$, see (Higgins, 1996). The classical Whittaker--Shannon--Kotel'nikov sampling theorem states that
any $f\in PW^2_{\omega}$ can be reconstructed at arbitrary point $x$ by its values at points $\frac {n\pi}w:$
   \begin{equation}\label{1}
      f(x) = \sum_{n=-\infty}^{\infty}f \left( \frac {n\pi}w\right) \frac{\sin (wx-n \pi)}{wx-n \pi}\,.
   \end{equation}
Although this formula yields an explicit reconstruction of $f(x)$, it is usually considered to be of theoretical interest only, because it
requires computations of infinite sums. In practice, we truncate the series in (\ref{1}). Typically, to estimate $f(x)$ we should only use the
values  of the function at points around $x,$ for example $\{ n: \left|n-\frac{\omega x}{\pi}\right|\le N\}.$  The sampling size $N$ is determined
by the relative error accepted in the reconstruction. Hence the reconstruction possesses time adapted sampling size. Thus the error analysis
plays a crucial role in setting up the interpolation formula.

The time shifted truncation version of (\ref{1})
   \[ Y_N(f;x):=\sum_{\left|n-\frac{\omega x}{\pi}\right|\le N}f \left( \frac {n\pi}w\right) \frac{\sin (wx-n \pi)}{wx-n \pi}\, \]
plays an important role in applications, when the appropriate remainder is conveniently bounded by a term which vanishes with the growing
sampling size parameter $ N.$ The scheme is called  a time shifted sampling, since it depends on the location $x$ of the reconstruction. For a
study of time shifted sampling schemes and numerous references, see (Flornes {\em et al.}, 1999); (Micchelli {\em et al.}, 2009); (Olenko and
Pog\'any, 2006, 2007, 2010).

Then the second problem arises. By physical and applications reasons the measured samples in
practice may not be the values of the measured function $f$ precisely at the sample time $\frac {n\pi}w$, but only the local
average of the function $f$ around $\frac {n\pi}w.$ So, the measured sample values are
   \[ \langle f, u_n\rangle_{U_n} = \int_{U_n} f(x)u_n(x){\rm d}x, \quad U_n = {\rm supp}(u_n)\subseteq \big[t_n - \sigma_n',
                 t_n+\sigma_n''\big]\, \]
for some sequence $\mathbf u :=\big( u_n(x)\big)_{n \in \mathbb Z}$ of non--negative, norma\-lized, that is $\langle 1, u_n\rangle \equiv 1$,
averaging functions.

The local averaging models and methods were introduced by Gr\"ochenig (1992), and developed by Butzer and Lei (1998, 2000). Recently Sun and
Zhou (2002, 2003), gave some results in this direction, while the stochastic counterpart of the average sampling was intensively
studied by He {\em et al.} (2006, 2007), see also their references inside. The interested reader can consult, for example, the exhaustive list
of references in (He {\em at al.}, 2007) too.

However, the mentioned stochastic average sampling results were restricted to weakly statio\-nary stochastic processes $\xi,$ while the
approximation ave\-rage sampling sums were used around the origin. The obtained error reconstruction bounds were derived under various
decay assumptions on the function $f$ for the deterministic case or covariance functions for the stochastic case. The error
reconstruction bounds were not sharp in the sense that they are constant and do not depend on the time variable $x.$ It has to be
mentioned as well that no {\em extremal functions} were given in the listed works for the average approximation problem.

The third problem is that the stationarity assumption is unacceptable for many practical problems in signal processing. In recent years there has
been growing interest in various specific models of nonstationarity. The concept of harmonizability is a natural generalization of stationarity
that includes a large class of nonstationary processes, while retaining some spectral representations of stationary processes, see
section~\ref{sec2}.  A wide class of real systems has the output which is a harmonizable process, see (Bochner, 1956). Our intention is to
extend the mentioned above stochastic sampling results to the time shifted average sampling, considered for harmonizable processes.
New techniques are required to obtain the desired results.

The paper addresses these three problems. We study the time shifted finite average sampling sums
   \[ {\mathcal A}_{\mathbf u, N}(\xi;t):=\sum_{\left|n-\frac{\omega t}{\pi}\right|\le N}
             \langle \xi, u_n\rangle_{\mathfrak J_n(t)}\cdot
             \frac{\sin (wt-n \pi)}{wt-n \pi}\, ,\]
in approximating the initial stochastic signal $\xi(t)$  by the weighted averages of $\xi$ over $ \mathfrak J_n(t) :=
\big[ n\pi/w-\sigma'_n(t), n\pi/w+\sigma''_n(t) \big],\, 0 \le \sigma'_n(t), \sigma''_n(t) \le \pi/(2w)$. The aim of the article
is to derive some reasonably simple
efficient truncation error upper bounds appearing in the approximation $\xi(t) \approx {\mathcal A}_{\mathbf u, N}(\xi;t)$. The mean--square
sampling truncation error
   \[ \mathfrak T_{\mathbf u, N}(\xi;t) := \mathsf E\big|\xi(t) - {\mathcal A}_{\mathbf u, N}(\xi;t)\big|^2\,\]
is investigated. We do not assume any decay conditions and construct sharp upper bounds on $\mathfrak T_{\mathbf u, N}(\xi;t)$ adapted to the
time variable $t.$

In comparison with all previous approaches our method is superior in obtaining very simple {\em minimal, time shifted} truncation error
upper bounds in interpolating  the harmonizable class process with Paley--Wiener class kernel function in its spectral representation.

The organization of this article is the following. In Section 2 we introduce the necessary background from the theory of
harmonizable stochastic processes. In the third section some Piranshvili's results are adapted to the time shifted sampling procedures.
In \S 4 we formulate and prove some auxiliary results on stochastic shifted average sampling. This section also contains the
main theorem and its corollaries. Conclusions are made in section~5.

\section{A Brief Review of Piranashvili Processes}\label{sec2}

In this section we give some key results concerning harmonizable processes and their spectral representations with respect to bimeasures.

Let $\{ \xi(t),  t \in \mathbb R\}$ be a centered second order random process defined on certain fixed probability space
$(\Omega, \mathfrak F, \mathsf P).$ Further, let the process $\xi(t)$ have  a covariance function (associated to some domain
$\Lambda \subseteq \mathbb R$ with some sigma--algebra $\sigma(\Lambda)$) in the form:
   \begin{equation} \label{I1}
      B(t,s) = \int_{\Lambda} \int_{\Lambda} f(t,\lambda)f^{*} (s, \mu) F_\xi({\rm d} \lambda, {\rm d} \mu),
   \end{equation}
where, for each $\lambda\in \Lambda,$ $f(\cdot,\lambda)$ can be extended to the complex plane as an complex analytic exponentially bounded kernel function, that is, for some $M>0,$ $\alpha\in \mathbb R$
   \[ |f(t, \lambda)| \le M{\rm e}^{\alpha|t|},\]
while $F_\xi(\cdot, \cdot)$ is a positive definite measure on $\mathbb R^2.$ The total variation $\|F_\xi\|(\Lambda, \Lambda)$ of the
spectral distribution function $F_\xi$ satisfies
   \[ \|F_\xi\|(\Lambda, \Lambda) = \int_\Lambda\int_\Lambda\big|F_\xi({\rm d}\lambda, {\rm d}\mu)\big|<\infty .\]
Notice that the sample function $\xi(t) \equiv \xi(t, \omega_0)$ and $f(t, \lambda)$ possess the same exponential types, see
(Belyaev, 1959; Theorem 4), and (Piranashvili, 1967; Theorem 3). Then, by the Karhunen--Cram\'er theorem the process $\xi(t)$ has the
spectral representation
   \begin{equation} \label{I2}
      \xi(t) = \int_\Lambda f(t, \lambda)Z_\xi({\rm d}\lambda)\, ,
   \end{equation}
where $Z_\xi(\cdot)$ is a stochastic measure and
   \[F_\xi(S_1,S_2) = {\mathsf E}Z_\xi(S_1)Z_\xi^*(S_2), \qquad S_1, S_2 \subseteq \sigma(\Lambda).\]
Such a process will be called Piranashvili process in the sequel, see (Piranashvili, 1967), and (Pog\'any, 1999).

Being $f(t, \lambda)$ entire, it possesses the Maclaurin expansion
   \[ f(t, \lambda) = \sum_{n=0}^{\infty} f^{(n)}(0, \lambda) t^{n}/n!\,.\]
Let
   \[ \gamma := \sup_\Lambda c(\lambda) =  \sup_\Lambda \overline{\lim_n} \sqrt[n]{|f^{(n)}(0, \lambda)|}<\infty\, .\]
Since the exponential type of $f(t, \lambda)$ is equal to $\gamma$, then for all $w> \gamma$ there holds
   \begin{equation} \label{I4}
      \xi (t) = \sum_{n \in {\mathbb Z}}^{} \xi \left( \frac {n\pi}w\right) \frac{\sin (wt-n \pi)}{wt-n \pi} ,
   \end{equation}
and the series converges uniformly in the mean square and almost surely {(Piranashvili,\! 1967; Theorem 1).} This result we call Whittaker--Kotel'\-nikov--Shannon
(WKS) stochastic sampling theorem (Pog\'any, 1999).

The class of Piranashvili processes includes various well known subclasses of stochastic processes. Some particular cases of Piranashvili
processes  are listen below. Specifying $F_\xi(x, y) = \delta_{xy}F_\xi(x)$ in (\ref{I1}) one easily concludes the Karhunen--representation of
the covariance function
   \[ B(t,s) = \int_{\Lambda} f(t, \lambda)f^{*}(s, \lambda)F_\xi({\rm d} \lambda) .\]
Also, putting $f(t, \lambda)= e^{it\lambda}$ in (\ref{I1}) one gets the Lo\`eve-representation:
   \[ B(t,s) = \int_{\Lambda} \int_{\Lambda} e^{i(t \lambda - s \mu)}F_\xi({\rm d} \lambda, {\rm d} \mu) .\]
Here is $c(\lambda) = |\lambda|$ and, therefore, WKS--formula (\ref{I4}) holds for all $w> \gamma = \sup|\Lambda|$.

Note that the Karhunen process with the Fourier kernel $f(t, \lambda)= e^{it\lambda}$ is the weakly stationary stochastic process
having the covariance
   \[ B(\tau) = \int_\Lambda e^{i\tau \lambda}F_\xi({\rm d}\lambda), \qquad \tau = t-s. \]
A deeper insight into different kinds of harmonizabilities is presented in (Kakihara, 1997); (Priestley, 1988); (Rao, 1982), and the
related references therein. Finally, using $\Lambda = [-w,w]$ for some finite $w$ in this consideration, we get the band--limited
variants of processes of the same kind.

\section{Truncation Errors in Sampling of Piranashvili\\ Processes}

In this section we adapt some Piranshvili's results on stochastic sampling to the time shifted
sampling procedures. First of all we introduce few auxiliary results. Let $N_x$ stand for
the integer nearest to $xw/\pi, x\in \mathbb R$, and let us denote
   \begin{align*}
               &\qquad \qquad \mathbb I_N(x) := \{ n\in \mathbb Z \colon |xw/\pi - n|\le N\}\,,\\
               & \Gamma_N(x) := \Big\{z \in \mathbb C\colon |z - N_x|< \big(N+\frac12\big)\frac \pi w\Big\}, \quad N\in \mathbb N.
   \end{align*}
In what follows,  the series
   \[ \lambda(q) := \sum_{n=1}^\infty \frac1{(2n-1)^q}\]
stands for the Dirichlet lambda function.

\begin{thm}\label{teor1} Let $f(z)$ be an entire  bounded on the real axis function of exponential type $\gamma<w.$ Denote
   \[ L_f := \sup_{\mathbb R} \big|f(x)\big|, \quad
         L_0(z) := \frac{4w L_f |\sin(wz)|}{\pi (w-\gamma)\big(1-e^{-\pi}\big)}\, .\]
Then for all $z\in \Gamma_N(x)$ and large enough $N\in \mathbb
N$ it holds
   \begin{equation} \label{pir1}
      \Big| \sum_{\mathbb Z\setminus \mathbb I_N(x)} f\left( \frac{n\pi}{w}\right)\cdot
                   \frac{\sin (wz-n \pi)}{wz-n \pi}\Big|  < \frac{L_0(z)}N\, .
   \end{equation}
\end{thm}

\noindent{\em Proof.} The estimate can be obtained by the contour integration method.  We adapt Piranashvili's results
see (Piranashvili, 1967; p. 648) to the set $\mathbb I_N(x)$ in the following way. We use the notation $\partial\Gamma_N(x)$ for the boundary
of the set  $\Gamma_N(x).$ By Cauchy's residue theorem for the contour integral
   \[ \frac1{2\pi {\rm i}}\oint_{\partial\Gamma_N(x)} \frac{f(\zeta)}{\sin(w\zeta)}\, \frac{{\rm d}\zeta}{\zeta-x}\, ,\]
we obtain the left side of the inequality~(\ref{pir1}). Then we derive the right hand side upper bound by using the classical upper estimate
for entire functions of exponential type $\gamma<w$, see (Achieser, 1992).
\hfill $\Box$

Here and in what follows we denote by
   \begin{equation} \label{trunc}
      Y_N(\xi;t) := \sum_{\mathbb I_N(t)}\xi \left( \frac{n\pi}w\right) \frac{\sin (wt-n \pi)}{wt-n \pi}
   \end{equation}
the time shifted truncated WKS restoration sum for the stochastic process $\xi.$

\begin{thm}\label{teor2} Let $\xi(t)$ be a Piranashvili process with exponentially bounded kernel function $f(t,\lambda)$ and
let
   \begin{align} \label{II21}
      \widetilde L_f := \sup_{\mathbb R} \sup_{\Lambda}|f(t,\lambda)|, \qquad 
      \widetilde L_0(t) := \frac{4 \widetilde L_f w\,|\sin(wt)|}{\pi (w-\gamma)\big(1-e^{-\pi}\big)}\, .
   \end{align}
Then for all $t \in \mathbb R$ and large enough $N,$ we have
   \begin{equation} \label{II3}
      \mathsf E \big| \xi(t) - Y_N(\xi; t)\big|^2 < \frac{\widetilde L_0^2(t)}{N^2}\,
                \|F_\xi\|(\Lambda, \Lambda)\, .
   \end{equation}
\end{thm}

\noindent{\em Proof.} If $\xi(t)$ is a Piranashvili process, then by  the spectral representation formula (\ref{I2}) we obtain that
the truncated sampling sum (\ref{trunc}) can be expressed (in the mean square sense) in the form
   \[ Y_N(\xi; t) = \int_\Lambda \sum_{\mathbb I_N(t)} f \left( \frac{n\pi}w, \lambda \right)
                    \frac{\sin (wt-n \pi)}{wt-n \pi}\, Z_\xi({\rm d}\lambda)
                  = \int_\Lambda  Y_N(f;t)\,Z_\xi({\rm d}\lambda) \, ,\]
where
   \[  Y_N(f;t) = \sum_{\mathbb I_N(t)}f \left( \frac{n\pi}w; \lambda\right) \frac{\sin (wt-n \pi)}{wt-n \pi}\, .\]
Now, combining the above formul{\ae} with (\ref{I2}) yields that
   \begin{align} \label{pir2}
      &\mathsf E \big| \xi(t) - Y_N(\xi; t)\big|^2 =
            \mathsf E \Big| \int_\Lambda \Big( f(t, \lambda) - Y_N(f;t)\Big)\, Z_\xi({\rm d}\lambda)\Big|^2 \nonumber \\
         &= \int_\Lambda \int_\Lambda \Big( f(t, \lambda) - Y_N(f;t)\Big) \cdot
            \Big( f^*(t, \mu) - Y_N^*(f;t)\Big)\, \mathsf EZ_\xi({\rm d}\lambda)Z_\xi^*({\rm d}\mu) \nonumber \\
         &= \int_\Lambda \int_\Lambda \Big( f(t, \lambda) - Y_N(f;t)\Big) \cdot
            \Big( f^*(t, \mu) - Y_N^*(f;t)\Big)\, F_\xi({\rm d}\lambda, {\rm d}\mu)\, .
    \end{align}
Denote
   \[ s_N(t) = \sup_{\lambda \in \Lambda} \big|f(t, \lambda) - Y_N(f;t)\big|\, .\]
Applying Theorem~{\rm \ref{teor1}}
twice to (\ref{pir2}), we deduce that
   \[ \mathsf E \big| \xi(t) - Y_N(\xi; t)\big|^2 \le s_N^2(t) \int_\lambda\int_\lambda \big|F_\xi({\rm d}\lambda, {\rm d}\mu)\big|
      \le \left( \frac{\widetilde L_0(t)}N \right)^2\, \int_\lambda\int_\lambda \big|F_\xi({\rm d}\lambda, {\rm d}\mu)\big|\, ,\]
completing the proof.
\hfill $\Box$

\begin{rmk}\label{rem1} Let us point out that the straightforward consequence of {\rm(\ref{II3})} is not only the exact $L_2$--restoration of
the initial Piranashvili--type harmonizable process $\xi$ by the sequence of approximants $Y_N(\xi; t)$ when $N \to \infty$, but
the {\rm a.s.} reconstruction as well. It follows immediately from
   \[ \mathsf E \big| \xi(t) - Y_N(\xi; t)\big|^2 = \mathcal O\big(N^{-2}\big)\, ,\]
and the  Borel--Cantelli lemma.
\end{rmk}

\section{Time Shifted Average Sampling}

In this section we derive some simple efficient mean square truncation error upper bounds appearing in the
approximation $\xi(t) \approx {\mathcal A}_{\mathbf u, N}(\xi;t)$.

Instead of the approach used by Song {\em et al.} (2007), we take time shifted  finite average sampling sums (consult
(Olenko and Pog\'any, 2006, 2007)) in approximating the initial stochastic signal $\xi$. First, we consider the weighted averages of
$\xi(t)$ over $ \mathfrak J_n(t)= \big[ n\pi/w-\sigma'_n(t), n\pi/w+\sigma''_n(t) \big],\, 0 \le \sigma'_n(t), \sigma''_n(t) \le \pi/(2w)$ instead of the measured values at $n\pi/w,\, n\in \mathbb I_N(t)\,.$

Suppose that the following assumption holds true:
   \[\sigma := \sup_{\mathbb R}\max_{n\in\mathbb I_N(t)}\max \big(\sigma'_n(t),\sigma''_n(t)\big)<\infty\, .\]
Let us define the time shifted average sampling approximation sum in the form
   \[{\mathcal A}_{\mathbf u}(\xi;t) = \sum_{\mathbb Z}
             \langle \xi, u_n\rangle_{\mathfrak J_n(t)} \cdot \frac{\sin (wt-n \pi)}{wt-n \pi}\, ,\]
and its truncated variant as
   \[ {\mathcal A}_{\mathbf u, N}(\xi;t) = \sum_{\mathbb I_N(t)}
             \langle \xi, u_n\rangle_{\mathfrak J_n(t)}\cdot
             \frac{\sin (wt-n \pi)}{wt-n \pi}\, .\]
We will study the mean--square, time shifted average sampling truncation error
   \[\mathfrak T_{\mathbf u, N}(\xi;t) = \mathsf E\big|\xi(t) - {\mathcal A}_{\mathbf u, N}(\xi;t)\big|^2\,.\]

\begin{lem}\label{lem1} Let $\xi(t)$ be a Piranashvili process with the covariance $B(t,s)\in C^2(\mathbb R)$, satisfying
$\sup_{\mathbb R}\big|B''(t,t)\big|<\infty$. Let $(p,q)$ be a conjugated H\"older pair of exponents:
   \[ \frac1p+\frac1q = 1,\qquad p >1\,.\]
Then
   \begin{equation} \label{II5}
      \mathsf E \big|Y_N(\xi;t) - \mathcal A_{\mathbf u, N}(\xi;t)\big|^2
                 \le \sigma^2\,  {\mathsf C_q(t)} \sup_{\mathbb R}\big|B''(t,t)\big|
                       \cdot (2N+1)^{2/p}\, ,
   \end{equation}
where
   \begin{equation} \label{II6}
      \mathsf C_q(t):= \left(1+ \frac{2^{q+1}\,|\sin(wt)|^q}{\pi^q}\, \lambda(q)\right)^{2/q}.
   \end{equation}
\end{lem}

\noindent{\em Proof.}
Let us note that the first order difference
$\Delta_{x,y}B$ see (Habib and  Cambanis, 1981) of $B(t,s)$ on the plane satisfies
   \begin{align} \label{II4}
      \big( \Delta_{x,y}B\big)(t,s) &= B(t+x, s+y) - B(t+x, s) - B(t,s+y) + B(t,s) \nonumber \\
            & = \int_0^x\int_0^y \frac{\partial^2}{\partial u\partial v}B\big(t+u,s+v\big)\, {\rm d}v{\rm d}u\, .
   \end{align}
Having in mind (\ref{I1}), the properties of the sequence $\mathbf u$ of averaging functions and (\ref{II4}), we obtain
   \begin{align*}
      & \hspace{4pt} \mathsf E \big|Y_N(\xi;t) - \mathcal A_{\mathbf u, N}(\xi;t)\big|^2 = \mathsf E \Big| \sum_{n\in\mathbb I_N(t)}
                             \langle \xi\big( \frac{n\pi}w\big) - \xi, u_n\rangle_{\mathfrak J_n(t)}
                       \cdot \frac{\sin (wt-n \pi)}{wt-n \pi}\Big|^2  \\
                    &\quad = \sum_{n,m\in\mathbb  I_N^2(t)}  \frac{\sin (wt-n \pi)}{wt-n \pi}\,
                            \frac{\sin (wt-m \pi)}{wt-m \pi} \int_{-\sigma'_n(t)}^{\sigma''_n(t)}
                             \int_{-\sigma'_m(t)}^{\sigma''_m(t)} u_n(x+n\frac\pi w)u_m(y+m\frac\pi w) \\
                             & \quad \times  \int_0^x\int_0^y\, \frac{\partial^2}{\partial u\partial v}
                            B\big(u+n\frac\pi w,v+m\frac\pi v\big){\rm d}v{\rm d}u {\rm d}x{\rm d}y  \le \sum_{\mathbb I_N^2(t)}
                            \Big|\frac{\sin (wt-n \pi)}{wt-n \pi}\Big|\,
                            \Big| \frac{\sin (wt-m \pi)}{wt-m \pi} \Big|  \\
                    & \quad \times \sup_{|x|,|y|\le \sigma}\Bigg| \int_0^x\int_0^y\, \frac{\partial^2}{\partial u\partial v}
                            B\big(u+n\frac\pi w,v+m\frac\pi v\big){\rm d}v{\rm d}u \Bigg|
   \end{align*}
being $\mathbf u$ normalized. For the sake of brevity, let us denote by $H_\sigma(n,m)$ the {\rm sup}--term in the last display. Then,
by the H\"older inequality with conjugate exponents $p>1$ and $q$, we get
   \[ \mathsf E \big|Y_N(\xi;t) - \mathcal A_{\mathbf u, N}(\xi;t)\big|^2
                    \le \Bigg\{\sum_{n,m\in \mathbb I_N^2(t)} H_\sigma^p(n,m)\Bigg\}^{1/p}
                       \Bigg\{\sum_{n\in\mathbb I_N(t)} \Big|\frac{\sin (wt-n \pi)}{wt-n \pi}\Big|^q\Bigg\}^{2/q}\, . \]
It is not hard to see that for all $n,m\in \mathbb I_N(t)$ there holds
   \begin{equation} \label{sig}
      H_\sigma(n,m) \le \sigma^2 \sup_{\mathbb R^2} \Big| \frac{\partial^2 B(t,s)}{\partial t \partial s}\Big| \, .
   \end{equation}
Applying the Cauchy--Bunyakovsky--Schwarz inequality to $\partial^2B$, we deduce
   \[ \sup_{\mathbb R^2} \Big| \frac{\partial^2 B(t,s)}{\partial t \partial s}\Big| \le
                               \sup_{\mathbb R} \Big| \frac{\partial^2 B(t,t)}{\partial t^2}\Big|
                            = \sup_{\mathbb R}|B''(t,t)|\, .\]
It remains to evaluate the sum of the $q$th power of the {\rm sinc}--functions. Since
   \[ \left|\frac{\sin(wt-N_t\pi)}{wt-N_t\pi}\right| \le 1\]
we conclude
   \begin{align*}
      \sum_{n\in\mathbb I_N(t)} \Big|\frac{\sin (wt-n \pi)}{wt-n \pi}\Big|^q &\le 1 + C(t)\cdot\sup_{\delta\in[0,\frac{1}{2}]}\,\sum_{n=1}^N
              \Big\{ \frac1{(n-\delta)^q}+\frac1{(n+\delta)^q}\Big\} \\
         &<   1+2C(t)\cdot\sum_{n=1}^\infty \frac1{(n-1/2)^q} =  1+ 2^{q+1}\lambda(q)\,C(t)\, ,
   \end{align*}
where
   \[ C(t): = \frac{|\sin(wt)|^q}{\pi^q}\,.\]
Collecting all these estimates, we deduce (\ref{II5}).
\hfill $\Box$

We are ready to formulate our main upper bound result for the mean square, time shifted average sampling truncation error
$\mathfrak T_{\mathbf u, N}(\xi;t)$. The almost sure sense restoration procedure will be treated too.

Since we use the average sampling sum $\mathcal A_{\mathbf u, N}(\xi;t)$ instead of $Y_N(\xi;t),$ to obtain asymptotically vanishing
$\mathfrak T_{\mathbf u, N}(\xi;t),$ it is not enough letting $N \to \infty$ as in Remark~\ref{rem1}. For the average sampling reconstruction
we need some additional conditions upon $w$ or $\sigma$ to guarantee smaller average intervals for larger/denser sampling grids.

\begin{thm}\label{teor3} Let the assumptions of  Theorem {\rm \ref{teor2}} and Lemma {\rm \ref{lem1}} hold true. Then
   \[ \mathfrak T_{\mathbf u, N}(\xi;t) \le \frac{2\widetilde L_0^2(t)}{N^2}\, \|F_\xi\|(\Lambda, \Lambda)
             + 2\sigma^2\,\mathsf C_q(t)\sup_{\mathbb R}\big|B''(t,t)\big| \cdot (2N+1)^{2/p}\, ,\]
where $\widetilde L_0$ and $\mathsf C_q(t)$ are given respectively by {\rm (\ref{II21})} and {\rm (\ref{II6})}.

If $\sigma =  o\big(N^{-1/p}\big),$ then
   \[ \lim_{N\to \infty} {\mathcal A}_{\mathbf u, N}(\xi;t) = \xi(t)\ \ \mbox{in mean square}.\]

Moreover, if $\sigma = \mathcal O\big(N^{-1/2-1/p- \varepsilon}\big),\, \varepsilon>0$, then
   \begin{equation} \label{II8}
      \mathsf P\big\{ \lim_{N\to \infty} {\mathcal A}_{\mathbf u, N}(\xi;t) = \xi(t)\big\} = 1
   \end{equation}
for all $t\in \mathbb R$.
\end{thm}

\noindent{\em Proof.} By direct calculation we deduce
   \begin{align*}
      \mathfrak T_{\mathbf u, N}(\xi;t)& = \mathsf E\big|\xi(t) - {\mathcal A}_{\mathbf u, N}(\xi;t)\big|^2
      =   \mathsf E\big|\xi(t) - Y_N(\xi;t)+ Y_N(\xi;t) - {\mathcal A}_{\mathbf u, N}(\xi;t)\big|^2 \\
          &\le 2 \mathsf E\big|\xi(t) - Y_N(\xi;t)\big|^2 + 2\mathsf E\big| Y_N(\xi;t) - {\mathcal A}_{\mathbf u, N}(\xi;t)\big|^2\, .
   \end{align*}
Thus, we get the asserted upper bound by (\ref{II3}) and (\ref{II5}). Hence, for $\sigma =  o\big(N^{-1/p}\big),$ we obtain
convergence in mean square.

To derive (\ref{II8}), we apply the Chebyshov inequality:
   \[ \mathsf P_N := \mathsf P\big\{ \big|\xi(t) - {\mathcal A}_{\mathbf u, N}(\xi;t)\big|\ge \eta\big\} \le
      \eta^{-2} \mathfrak T_{\mathbf u, N}(\xi;t)\, .\]
Since $\widetilde L_0(t) = \mathcal O(1)$ as $N \to \infty$, we obtain
   \[ \sum_{\mathbb N} \mathsf P_N \le K\,\sum_{\mathbb N}\left( \frac1{N^2} + \sigma^2\,(2N+1)^{2/p}\right) <\infty,\]
where $K$ is a suitable absolute constant.

An application of the Borel--Cantelli lemma results in the a.s. convergence (\ref{II8}), which completes the proof.
\hfill $\Box$

\begin{rmk}
{\rm Theorem~{\rm\ref{teor3}}} gives new truncation error upper bounds in time shifted average sampling restorations for wide classes of
harmonizable processes. The estimate~{\rm(\ref{sig})} and the upper bound for $\sum_{n,m\in \mathbb I_N^2(t)} H_\sigma^p(n,m)$ in
{\rm Lemma~{\rm\ref{lem1}}} can be specified for particular classes of stochastic processes which correlation or kernel functions decay
rates are known. Thus, one can use the suggested approach and obtained upper bounds to sharpen the results, mentioned in the
introduction.
\end{rmk}

\begin{rmk}
The upper bounds in {\rm Lemma~{\rm\ref{lem1}}} and {\rm Theorem~{\rm\ref{teor3}}} are sharp. It is easy to check choosing $t=\frac{n\pi}{\omega}$
and trivial stochastic process $\xi(t)\equiv\xi_0,$ where $\xi_0$ is a random variable with finite variance.
\end{rmk}

\begin{rmk}  By obvious reasons, in many applications they restrict the study to $\sigma \le \frac \pi{2w}.$ In this case we can replace
the estimate~{\rm(\ref{sig})} by
   \[ H_\sigma(n,m) \le \frac{\pi^2}{4w^2}\, \sup_{\mathbb R^2} \Big| \frac{\partial^2 B(t,s)}{\partial t \partial s}\Big|\, .\]
Hence, instead of the  estimate~{\rm(\ref{II5})} in {\rm Lemma~{\rm\ref{lem1}}} we obtain
   \[\mathsf E \big|Y_N(\xi;t) - \mathcal A_{\mathbf u, N}(\xi;t)\big|^2
                 \le \frac{\mathsf \pi^2\,C_q(t)}{4w^2}\sup_{\mathbb R}\big|B''(t,t)\big| \cdot (2N+1)^{2/p}\,. \]
In this case {\rm Theorem {\rm \ref{teor3}}} ensures the perfect time shifted average sampling restoration in the mean square sense when
$w = \mathcal O\big(N^{1/p+\varepsilon}\big),\, \varepsilon>0.$ The a.s. sense restoration~{\rm(\ref{II8})} requires a stronger assumption.
It holds when $w = \mathcal O\big(N^{1/2+1/p+ \varepsilon}\big)$.

In this case {\rm Theorem {\rm \ref{teor3}}} gives the rate of mean square convergence when the observations get dense in the whole
region (infill asymptotics).
\end{rmk}

\begin{rmk} In both cases we use the so-called {\em approximate sampling procedure}, when in the restoration procedure $\sigma\to 0$ or
$w \to \infty$ in some fashion. The consequence of this is that we have to restrict ourselves to the case $\Lambda = \mathbb R.$ Therefore
we deal with the non--bandlimited Piranashvili type harmonizable process case.

The importance of the approximate sampling procedures for investigations of aliasing errors in sampling restorations and various conditions
on the joint asymptotic behaviour of $N$ and $w$ were discussed in detail in {\rm (Olenko and  Pog\'any, 2006)}.
\end{rmk}

\section{Conclusions}
We have analyzed truncation error upper bounds in time shifted average sampling restorations for the stochastic initial signal
case. The general class of harmonizable stochastic processes was studied. The convergence of the truncation error to zero was discussed.
The results are obtained under simple conditions.  The conditions are weaker than those in the former literature: there are no any
assumptions on decay rates of correlation or kernel functions. The analysis is new and provides a constructive algorithm for
determining the number of terms in the sampling expansions to ensure the approximation of stochastic processes with given accuracy.

However, certain new questions immediately arise:
   \begin{enumerate}
      \item[(i)] to apply the developed techniques to $L_p$--processes using recent deterministic findings in (Olenko and Pog\'any, 2010, 2011);
      \item[(ii)] to obtain similar results for irregular/nonuniform sampling restorations using the methods exposed
      in (Olenko and  Pog\'any, 2003, 2011).
   \end{enumerate}

\section{Acknowledgements}
This work was partly supported by La Trobe University Research Grant--501821 "Sampling, wavelets and optimal stochastic modelling".

\vskip 3mm

\noindent \textbf{\Large References}
\vskip 3mm

\noindent  Achieser, N. I. (1992). {\it Theory of Approximation.} New York: Dover Publications, Inc.

\noindent  Belyaev, Yu.K. (1959). Analytical random processes. {\it  Theory Probab. Appl.}  IV(4):402--409.

\noindent   Bochner, S. (1956). Stationarity, boundedness, almost periodicity of random-valued functions. In: {\it Proc. Third Berkeley  Symp. on Math. Statist. and Prob.}  Vol. 2. Univ. of Calif. Press., 7--27.

\noindent  Butzer, P.L.,  Lei, J. (1998). Errors in truncated sampling series with measured sampled values for non-necessarily
bandlimited functions.~{\it Funct.~Approx.~Comment.~Math.}~26:25--39.

\noindent  Butzer, P.L.,  Lei, J. (2000).  Approximation of signals using measured sampled values and error analysis. {\it Commun.
Appl. Anal.}  4:245--255.

\noindent   Gr\"ochenig, K. (1992). Reconstruction algorithms in irregular sampling. {\it Math. Comp.}  59:181--194.

\noindent  Flornes, K. M.,  Lyubarskii, Yu., Seip, K. (1999). A direct interpolation method for irregular sampling. {\it Appl. Comput. Harmon. Anal.}  7(3):305--314.

\noindent  Habib, M.K.,  Cambanis, S. (1981). Sampling approximation for non--band--limited harmonizable random signals. {\em Inform. Sci.}  23:143--152.

\noindent  He, G.,  Song, Zh.,  Yang, D.,  Zhu, J. (2007). Truncation error estimate on random signals by local average. In: Shi, Y., et al., ed., {\it ICCS 2007, Part II, Lecture Notes in Computer Sciences 4488.}
Berlin: Springer-Verlag,  1075--1082.

\noindent   Higgins, J.R. (1996). {\it Sampling in Fourier and Signal Analysis: Foundations.} Oxford: Clarendon Press.

\noindent   Kakihara, Y. (1997). {\it Multidimensional Second Order Stochastic Processes.} Singapore: World Scientific.

\noindent  Micchelli, C. A.,  Xu, Yu., Zhang, H. (2009). Optimal learning of bandlimited functions from localized sampling. {\it J. Complexity.} 25(2):85--114.

\noindent  Olenko, A.,  Pog\'any, T. (2003). Direct Lagrange--Yen type interpolation of random fields.
{\it Theory Stoch. Process.}  9(25)(3-4):242--254.

\noindent   Olenko, A.,  Pog\'any, T. (2006). Time shifted aliasing error upper bounds for truncated sampling cardinal series. {\it J. Math. Anal. Appl.}  324(1):262--280.

\noindent   Olenko, A.,  Pog\'any, T. (2007). On sharp bounds for remainders in multidimensional sampling theorem. {\it Sampl. Theory Signal Image Process.} 6(3):249--272.

\noindent   Olenko, A.,  Pog\'any, T. (2010). Universal truncation error upper bounds in sampling restoration. {\it Georgian Math. J.}
17(4):765--786.

\noindent   Olenko, A.,  Pog\'any, T. (2011). Universal truncation error upper bounds in irregular sampling restoration.
(to appear in {\it Appl. Anal.})

\noindent  Piranashvili, Z. (1967). On the problem of interpolation of random processes. {\it Theory Probab. Appl.}
 XII(4):647--657.

\noindent  Pog\'any,  T. (1999). Almost sure sampling restoration of bandlimited stochastic signals.
In: Higgins J.R.,  Stens, R.L., ed., {\it Sampling Theory in Fourier and Signal Analysis: Advanced Topics.} Oxford
University Press, 203--232, 284--286.

\noindent   Priestley, M. (1988). {\it Non--linear and non--stationary time series.}~New~York:~Academic~Press.

\noindent   Rao,\! M.\! (1982).\! Harmonizable processes:~structure~theory.~{\it Enseign.~Math.(2)}~\mbox{28(3-4):295-351.}

\noindent   Song, Zh.,  Zhu, Z.,  He, G. (2006). Error estimate on non--bandlimited random signals
by local averages. In: Aleksandrov, V.N., et al., ed., {\it ICCS 2006, Part I,  Lecture Notes in
Computer Sciences 3991.} Berlin: Springer-Verlag, 822--825.

\noindent   Song, Zh.,  Yang, Sh.,  Zhou, X. (2006). Approximation of signals from local
averages. {\it Appl. Math. Lett.} 19:1414--1420.

\noindent   Song, Zh.,  Sun, W., Yang, Sh., Zhu, G. (2007). Approximation of weak sense
stationary stochastic processes from local averages. {\it Sci. China Ser. A.}  50(4):457--463.

\noindent  Sun, W.,  Zhou, X. (2002). Reconstruction of bandlimited signals from local averages. {\it IEEE Trans. Inform. Theory.}
48:2955--2963.

\noindent  Sun, W.,  Zhou, X. (2003). Reconstruction of functions in spline subspaces from local averages. {\it Proc. Amer.
Math. Soc.}  131:2561--2571.

\end{document}